\documentclass[a4paper, 11pt, oneside]{amsart}

\usepackage{a4,graphicx,mathrsfs,amssymb, tikz, cancel, soul, color}
\usepackage[colorlinks=true,breaklinks=true]{hyperref}
\usepackage{enumerate}
\usepackage{multirow}

\newtheorem{theorem}{Theorem}[section]

\newtheorem{proposition}[theorem]{Proposition}
\theoremstyle{definition}
\newtheorem{remark}[theorem]{Remark}
\newtheorem{definition}[theorem]{Definition}
\theoremstyle{remark}

\numberwithin{equation}{section}

\newcommand{\field}[1]{\mathbb{#1}}
\newcommand{\C}{\field{C}}

\newcommand{\N}{\field{N}}

\setstcolor{blue}

\allowdisplaybreaks

\begin{document}
\title[On  Trace Theorems and Poincare inequality]{On  Trace Theorems and Poincare inequality for one-dimensional Sobolev spaces}

\author{Bienvenido Barraza Mart\'inez}
\address{B.\ Barraza Mart\'inez, Universidad del Norte, Departamento de Matem\'aticas y Estad\'istica, Barranquilla, Colombia}
\email{bbarraza@uninorte.edu.co}

\author{Jonathan Gonz\'alez Ospino}
\address{J. Gonz\'alez Ospino, Universidad del Norte, Departamento de Matem\'aticas y Estad\'istica, Barranquilla, Colombia}
\email{gjonathan@uninorte.edu.co}

\author{Jairo Hern\'andez Monz\'on}
\address{J.\ Hern\'andez Monz\'on, Universidad del Norte, Departamento de Matem\'aticas y Estad\'istica, Barranquilla, Colombia}
\email{jahernan@uninorte.edu.co}

\renewcommand{\shortauthors}{B. Barraza, J. Gonz\'alez and J. Hern\'andez}


\maketitle

\begin{abstract}
In these notes, we present versions of trace theorems for Sobolev spaces over an interval in the real line, and also a one-dimensional version of the well-known Poincare inequality. \\
\end{abstract}

\section{Introduction}
The aim of these notes is to provide particular results relative to traces of functions in Sobolev spaces as well as a generalized version of Poincar\'e's inequality, all in one dimension, that can be cited in academic works where such results are needed for the analysis of transmission problems that model the dynamics of coupled structures in one dimension. Of course, such results already exist in the classical literature on Sobolev spaces (see e.g. \cite{Adams, Necas}), but mostly in multidimensional form, and we believe that it may be practical to have such concrete one-dimensional results at hand.\\

Throgout this notes let $a$ and $b$ be real numbers with $a<b$ and $\lambda:=b-a$. For $m\in\N$ and $1\leq p<\infty$ let $H^{m,p}(a,b)$ be defined as the completion of $C^m([a,b])$ with respect to the norm
$$ \Vert u\Vert_{m,p}:=\Big(\sum\limits_{j=0}^m\int_a^b \big\vert u^{(j)}(x)\big\vert^p\,dx\Big)^{1/p},$$
where $C^m([a,b])$ denotes the set of all complex valued $m$ times continuously differentiable fuctions in $[a,b]$, and $u^{(j)}$ is the derivative of order $j$ of $u$.\\

The elements belonging to $H^{m,p}(a,b)$ are indentified with functions in $L^p(a,b)$ whose derivatives up to order $m$ also belong to $L^p(a,b)$ (see \cite[Chapter 3]{Adams}).\\

As usual, if $p=2$ we write $H^m(a,b)$ instead $H^{m,2}(a,b)$ and we will write also $\Vert\cdot\Vert_{H^m(a,b)}$ instead $\Vert\cdot\Vert_{m,2}$.

 \section{Trace theorems for Sobolev spaces in one dimension}
 Next we will present versions, in the one dimensional case, of the well known trace theorems for Sobolev spaces $H^m(a,b)$, with $m\in\N$, endowed with its usual norm. For this, let us at first make precise some concepts.
 \begin{definition}
  We define the space $L^2(\{a,b\})$ as the space of all functions $\varphi:\{a,b\}\to \C$  endowed with the norm given by
  $$\Vert \varphi\Vert_{L^2(\{a,b\})} := \Big( \vert \varphi(a)\vert^2 + \vert \varphi(b)\vert^2\Big)^{1/2}.$$
 \end{definition}
\begin{proposition}
 $\big(L^2(\{a,b\}), \Vert \cdot\Vert_{L^2(\{a,b\})}\big)$ is a Banach space.
\end{proposition}
\begin{proof}
 Let $(\varphi_n)_{n\in\N}$ be a Cauchy sequence in $L^2(\{a,b\})$. Then $\big(\varphi_n(a)\big)_{n\in\N}$ and $\big(\varphi_n(b)\big)_{n\in\N}$ are Cauchy sequences in $\C$. Therefore, there exist $\alpha, \beta\in\C$ with $\varphi_n(a)\to \alpha$ and $\varphi_n(b)\to \beta$ in $\C$, when $n\to\infty$. Now we define $\varphi:\{a,b\}\to\C$ by $\varphi(a):=\alpha$ and $\varphi(b):=\beta$. Then we have $\Vert \varphi_n - \varphi \Vert_{L^2(\{a,b\})} \to 0$ when $n\to\infty$.
\end{proof}
Now we can establish the following trace theorem.
\begin{theorem}[Trace theorem]\label{Trace theoremo order zero}
The operator 
\begin{equation}\label{Eq_Trace_order0}
u \mapsto u\vert_{\{a,b\}} : C^1([a,b])\to L^2(\{a,b\})
\end{equation}
admits a unique extension to a bounded linear operator  
$$T_0:H^1(a,b)\to L^2(\{a,b\}).$$ 
\end{theorem}
\begin{proof}
 For $u\in C^1([a,b])$ and $a\leq \tau \leq b$ we have
 $$ u(a) = -\int_a^\tau u'(s)\,ds + u(\tau).$$
 Then, due to triangular and Cauchy-Schwarz inequalities, we obtain
 \begin{equation}\label{Eq1}
 \vert u(a) \vert \leq \int_a^b  \vert u'(s)\vert\,ds + \vert u(\tau)\vert \leq \sqrt{\lambda}\Big( \int_a^b \vert u'(s)\vert^2\,ds\Big)^{1/2} + \vert u(\tau)\vert.
 \end{equation}
 Therefore Young inequality yields
 \begin{equation*}
  \vert u(a) \vert^2  \leq 2\lambda\int_a^b \vert u'(s)\vert^2\,ds + 2\vert u(\tau)\vert^2, \qquad a\leq\tau\leq b.
 \end{equation*}
 Now, by integration with respect to $\tau$ in $[a,b]$ in the last inequality, we have
 \begin{equation*}
  \lambda\vert u(a) \vert^2 \leq 2\lambda^2\int_a^b \vert u'(s)\vert^2\,ds + 2\int_a^b\vert u(\tau)\vert^2\,d\tau
 \end{equation*}
and then
\begin{equation}\label{Eq2}
 \vert u(a) \vert^2 \leq 2\lambda\int_a^b \vert u'(s)\vert^2\,ds + \frac{2}{\lambda}\int_a^b\vert u(\tau)\vert^2\,d\tau.
\end{equation}
Noting that $u(b) = \int_\tau^b u'(s)\,ds + u(\tau)$ for $a\leq\tau\leq b$, we obtain similarly to \eqref{Eq1} that
$$\vert u(b) \vert  \leq \sqrt{\lambda}\Big( \int_a^b \vert u'(s)\vert^2\,ds\Big)^{1/2} + \vert u(\tau)\vert$$
and therefore we obtain again
\begin{equation}\label{Eq3}
 \vert u(b) \vert^2 \leq 2\lambda\int_a^b \vert u'(s)\vert^2\,ds + \frac{2}{\lambda}\int_a^b\vert u(\tau)\vert^2\,d\tau.
\end{equation}
Summing \eqref{Eq2} and \eqref{Eq3} we have
\begin{align}
\vert u(a) \vert^2 + \vert u(b) \vert^2 & \leq 4\lambda\int_a^b \vert u'(s)\vert^2\,ds + \frac{4}{\lambda}\int_a^b\vert u(\tau)\vert^2\,d\tau \nonumber\\
& \leq C_\lambda^2\Big( \int_a^b \vert u'(s)\vert^2\,ds + \int_a^b\vert u(\tau)\vert^2\,d\tau \Big)\label{Eq_induction1}
\end{align}
and therefore
\begin{equation}\label{Eq4}
\big\Vert u\vert_{\{a,b\}} \big\Vert_{L^2(\{a,b\})} \leq C_\lambda \Vert u \Vert_{H^1(a,b)}
\end{equation}
for all $u\in C^1([a,b])$, where, in this case\footnote{From now on, $C_\lambda$ will denotes several positive constants, which will only depend on $\lambda$.}, $C_\lambda := \Big( \max\Big\{4\lambda, \dfrac{4}{\lambda}\Big\} \Big)^{1/2}$.\\
Since $ C^1([a,b])$ is dense in $H^1(a,b)$ (indeed $H^1(a,b)$ is the completion of $ C^1([a,b])$ with respect to the usual norm in $H^1(a,b)$), there exists a unique linear extension $T_0:H^1(a,b)\to L^2(\{a,b\})$ of \eqref{Eq_Trace_order0} from $C^1([a,b])$ to $H^1(a,b)$, satisfying
\begin{equation}\label{Eq6}
\big\Vert T_0u \big\Vert_{L^2(\{a,b\})} \leq C_\lambda \Vert u \Vert_{H^1(a,b)} \qquad \big(u\in H^1(a,b)\big).
\end{equation}
 Then $T_0:H^1(a,b)\to L^2(\{a,b\})$ is a bounded linear operator and $T_0u = u\vert_{\{a,b\}}$ for all $u\in C^1([a,b])$.
\end{proof}
\begin{remark}
The operator $T_0:H^1(a,b)\to L^2(\{a,b\})$ is called \textit{Trace operator of order zero}. Usually this operator is also denoted $\gamma_0$.
\end{remark}
This last result can be generalized to traces of higher order.
\begin{theorem}[Higher order traces]
Let $m\in\N$, $m>1$. The operator 
\begin{equation}\label{Eq_trace_higher_1}
u \mapsto \big(u\vert_{\{a,b\}},u'\vert_{\{a,b\}},\cdots, u^{(m-1)}\vert_{\{a,b\}}\big): C^m([a,b])\to \prod\limits_{j=0}^{m-1}L^2(\{a,b\}) 
\end{equation}
admits a unique extension to a bounded linear  operator 
$$T_{m-1}:H^m(a,b)\to \prod\limits_{j=0}^{m-1}L^2(\{a,b\})$$ 
with $T_{m-1}u:=(\gamma_0u, \gamma_1u,\cdots,\gamma_{m-1}u)$, where  $\gamma_ju = u^{(j)}\vert_{\{a,b\}}$ for  $j=0,1,\dots,m-1$, and $u\in C^{m}([a,b])$. Here $\prod\limits_{j=0}^{m-1}L^2(\{a,b\})$ is the (cartesian) product of $m$ copies of $L^2(\{a,b\})$, endowed with the usual product topology. $T_{m-1}$ is called \emph{trace operator of order $m-1$}. 
\end{theorem}
\begin{proof}
Let $u\in C^m([a,b])$. Analogously to \eqref{Eq_induction1}, for $k\in \{0,1,\dots,m-1\}$ we have
\begin{equation*}
\vert u^{(k)}(a) \vert^2 + \vert u^{(k)}(b)\vert^2 \leq C_\lambda^2\Big( \int_a^b \vert u^{(k+1)}(t)\vert^2\,dt + \int_a^b\vert u^{(k)}(\tau)\vert^2\,d\tau \Big).
\end{equation*}
We recall that $C_\lambda$ denotes several constants, which depend only on $\lambda$. Then, adding these inequalities, we get
\begin{equation*}
\sum\limits_{k=0}^{m-1}\Big( \vert u^{(k)}(a) \vert^2 + \vert u^{(k)}(b)\vert^2\Big) \leq C_\lambda^2\sum\limits_{j=0}^{m}\int_a^b\vert u^{(j)}(t)\vert^2\,dt.
\end{equation*}
That is,
\begin{equation*}
\sum\limits_{k=0}^{m-1}\big\Vert u^{(k)}\vert_{\{a,b\}}\big\Vert_{L^2(\{a,b\})}^2 \leq C_\lambda^2\sum\limits_{j=0}^{m}\big\Vert u^{(j)}\big\Vert_{L^2(a,b)}^2.
\end{equation*}
Therefore,
\begin{align}
\Big\Vert \big(u\vert_{\{a,b\}},u'\vert_{\{a,b\}},\cdots, & u^{(m-1)}\vert_{\{a,b\}}\big) \Big\Vert_{\prod\limits_{j=0}^{m-1}L^2(\{a,b\})} \nonumber \\
& = \Big(\sum\limits_{k=0}^{m-1}\big\Vert u^{(k)}\vert_{\{a,b\}}\big\Vert_{L^2(\{a,b\})}^2\Big)^{1/2}\nonumber\\
& \leq C_\lambda \Big(\sum\limits_{j=0}^{m}\big\Vert u^{(j)}\big\Vert_{L^2(a,b)}^2\Big)^{1/2}\nonumber\\
& = C_\lambda \Vert u \Vert_{H^m(a,b)},\nonumber 
\end{align}
for all $u\in C^m([a,b])$.  The estimates above shows that 
\begin{equation}\label{Estimate_trace_Hm}
\big\Vert T_{m-1}u \big\Vert_{\prod\limits_{j=0}^{m-1}L^2(\{a,b\})} \leq C_\lambda \Vert u \Vert_{H^m(a,b)}\qquad (u\in C^m([a,b])). 
\end{equation}
Since $C^m([a,b])$ is dense in $H^m(a,b)$ ($H^m(a,b)$ is the completion of $C^m([a,b])$ with respect to the usual norm in $H^m(a,b)$), there exists a unique operator $T_{m-1}:H^m(a,b)\to \prod\limits_{j=0}^{m-1}L^2(\{a,b\})$ which extends \eqref{Eq_trace_higher_1} from   $C^m([a,b])$ to $H^m(a,b)$, satisfying
\begin{equation*}
\big\Vert T_{m-1}u \big\Vert_{\prod\limits_{j=0}^{m-1}L^2(\{a,b\})} \leq C_\lambda \Vert u \Vert_{H^m(a,b)}\qquad (u\in H^m(a,b)). 
\end{equation*}
If we define $T_{m-1}u:=(\gamma_0u, \gamma_1u,\cdots,\gamma_{m-1}u)$ for $u\in H^m(a,b)$, then we have, as extension of \eqref{Eq_trace_higher_1}, that
$\gamma_ju = u^{(j)}\vert_{\{a,b\}}$ for  $j=0,1,\dots,m-1$, and $u\in C^{m}([a,b])$.
\end{proof}
\begin{remark}
Since $m > (m-1)+\dfrac{1}{2}$, due to the Sobolev imbedding Theorem, we have $H^m(a,b)\hookrightarrow C^{m-1}([a,b])$. Therefore, for $u\in H^m(a,b)$ we can take its $C^{m-1}$ representative and then, the traces $u\vert_{\{a,b\}},u'\vert_{\{a,b\}},\cdots,u^{(m-1)}\vert_{\{a,b\}}$, exist in the classical sense. That is, it holds also
$$T_{m-1}u = \big(u\vert_{\{a,b\}},u'\vert_{\{a,b\}},\cdots,u^{(m-1)}\vert_{\{a,b\}}\big)$$ for all $u\in H^{m}(a,b)$.

\end{remark}
\section{Poincare inequality in one dimension}
In this section we show a generalized version of the well known Poincare inequality in the case of dimension one.
\begin{theorem}[Poincare inequality]
There exists a positive constant $C_P$ such that for $x_0\in\{a, b \}$ it holds
\begin{equation}\label{Eq_Poincare_inequaliy}
\Vert u \Vert_{L^2(a,b)} \leq C_P\Big( \Vert u'\Vert_{L^2(a,b)} + \vert T_0u(x_0)\vert \Big)
\end{equation}
for all $u\in H^1(a,b)$, where $T_0:H^1(a,b)\to L^2(\{a,b\})$ is the trace operator of order zero given in Theorem \ref{Trace theoremo order zero}.
\end{theorem}
\begin{proof}
Without lost of generality let us suppose that $x_0=a$ (the proof for $x_0=b$ is similar). Then, for $u\in C^1([a,b])$ and $a\leq x\leq b$, we have
\begin{equation*}
u(x) = \int_a^x u'(t)\,dt + u(a).
\end{equation*}
Therefore, applying triangular and Cauchy-Schwarz inequalities, we obtain
\begin{align*}
\vert u(x) \vert & \leq \int_a^x \vert u'(t) \vert\,dt + \vert u(a) \vert \\
& \leq \int_a^b \vert u'(t) \vert\,dt + \vert u(a) \vert \\
& \leq \sqrt{\lambda}\Big( \int_a^b \vert u'(t)\vert^2\,dt\Big)^{1/2} + \vert u(a) \vert.
\end{align*}
Then, Young inequality yields
\begin{equation*}
  \vert u(x) \vert^2  \leq 2\lambda\int_a^b \vert u'(t)\vert^2\,dt + 2\vert u(a)\vert^2, \qquad a\leq x\leq b.
 \end{equation*}
 We integrate both sides of the last inequality with respect to $x$ over the interval $[a,b]$ and obtain
 \begin{equation*}
 \int_a^b \vert u(x) \vert^2\,dx \leq 2\lambda^2\int_a^b \vert u'(t)\vert^2\,dt + 2\lambda\vert u(a)\vert^2.
 \end{equation*}
 That is, we have
 \begin{equation}\label{Eq_Poincare_smooth}
 \Vert u \Vert_{L^2(a,b)}^2 \leq 2\lambda^2\Vert u'\Vert_{L^2(a,b)}^2 + 2\lambda\vert u(a)\vert^2
 \end{equation}
 for all $u\in C^1([a,b])$.\\
 Now, let $u\in H^1(a,b)$ and $(u_n)_{n\in\N}$ a sequence of functions belonging to $C^1([a,b])$ such that $u_n\to u$ in $H^1(a,b)$ when $n\to \infty$. In virtue of \eqref{Eq_Poincare_smooth}, the inequality
 \begin{equation}\label{Eq_Poincare_smooth_sequence}
 \Vert u_n \Vert_{L^2(a,b)}^2 \leq 2\lambda^2\Vert u_n'\Vert_{L^2(a,b)}^2 + 2\lambda\vert u_n(a)\vert^2
 \end{equation}
 holds for all $n\in\N$. Due to the convergence of the sequence $(u_n)_{n\in\N}$  to $u$ in $H^1(a,b)$, it follows clearly that
 $$ \Vert u_n \Vert_{L^2(a,b)}^2 \to \Vert u \Vert_{L^2(a,b)}^2 \quad \text{and}\quad \Vert u_n'\Vert_{L^2(a,b)}^2 \to \Vert u'\Vert_{L^2(a,b)}^2 $$
 when $n\to \infty$. For the remaining term in \eqref{Eq_Poincare_smooth_sequence}, note that the Trace theorem (inequality \eqref{Eq6}) implies
 \begin{align*}
 \big\vert \vert u_n(a)\vert - \vert T_0u(a)\vert  \big\vert & \leq \vert T_0u_n(a) - T_0u(a) \vert \\
 & = \vert T_0(u_n - u)(a) \vert \\
 & \leq \Vert T_0(u_n-u)\Vert_{L^2(\{a,b\})}\\
 & \leq C \Vert u_n - u \Vert_{H^1(a,b)} \xrightarrow[n\to\infty]{} 0.
 \end{align*}
 Then, $ \vert u_n(a)\vert \xrightarrow[n\to\infty]{} \vert T_0u(a)\vert$.\\
 Taking $n\to\infty$ in \eqref{Eq_Poincare_smooth_sequence}, we obtain
 \begin{equation*}
 \Vert u \Vert_{L^2(a,b)}^2 \leq 2\lambda^2\Vert u'\Vert_{L^2(a,b)}^2 + 2\lambda\vert T_0u(a)\vert^2.
 \end{equation*}
 It follows that
 \begin{align*}
 \Vert u \Vert_{L^2(a,b)}^2 & \leq 2\lambda^2\Vert u'\Vert_{L^2(a,b)}^2 + 2\lambda\vert T_0u(a)\vert^2\\
 & \leq 2\max\{\lambda, \lambda^2\}\Big(\Vert u'\Vert_{L^2(a,b)}^2 + \vert T_0u(a)\vert^2 \Big)\\
 & \leq 2\max\{\lambda, \lambda^2\}\Big(\Vert u'\Vert_{L^2(a,b)} + \vert T_0u(a)\vert \Big)^2
 \end{align*}
 for all $u\in H^1(a,b)$. Taking square root in the last inequality, we obtain \eqref{Eq_Poincare_inequaliy} for $x_0=a$, with $C_P:=\sqrt{2}\big(\max\{\lambda, \lambda^2\}\big)^{1/2} = \sqrt{2}\max\{\sqrt{\lambda}, \lambda\} $. With similar calculations we obtain \eqref{Eq_Poincare_inequaliy} for $x_0=b$, with the same $C_P$.
\end{proof}
\begin{remark}
The Poincare inequality \eqref{Eq_Poincare_inequaliy} leads to the so-called \emph{Friedrichs} inequality (see \cite{Necas}, Theorem 1.9.) in dimension one:
\begin{equation}\label{Eq_Friedrichs inequality}
\Vert u\Vert_{H^1(a,b)} \leq \mathrm{const}\,\Big( \Vert u'\Vert_{L^2(a,b)}^2 + \vert T_0u(x_0)\vert^2 \Big)^{1/2}\qquad (u\in H^1(a,b)),
\end{equation}
which in turn implies the equivalence between the usual norm $\Vert\cdot\Vert_{H^1(a,b)}$ and the norm in $H^1(a,b)$ given by
$$u\mapsto \Big( \Vert u'\Vert_{L^2(a,b)}^2 + \vert T_0u(x_0)\vert^2 \Big)^{1/2}\qquad (u\in H^1(a,b)).$$
\end{remark}

\end{document}